\documentclass[preprint,12pt]{elsarticle}
\journal{ }

\usepackage{amsmath,amssymb,amsthm,a4wide}
\usepackage{latexsym}
\usepackage{indentfirst}
\usepackage{graphicx}
\usepackage{placeins}
\usepackage{enumerate}
\usepackage{booktabs}
\usepackage{algorithm}
\usepackage{algorithmic}
\usepackage{multirow}
\usepackage{subcaption}
\usepackage{optidef}
\usepackage{graphicx,color}
\usepackage{diagbox}
\usepackage[normalem]{ulem}

\newcommand{\R}{\mathbb{R}}

\begin{document}
\begin{frontmatter}
\title{Solving high-dimensional Fokker-Planck equation with functional hierarchical tensor}
\author[stanford]{Xun Tang}
\author[stanford]{Lexing Ying}

\affiliation[stanford]{organization={Stanford University},%Department and Organization
            % addressline={}, 
            city={Stanford},
            state={California},
            postcode={94305}, 
            country={USA}}

\begin{abstract}
    This work is concerned with solving high-dimensional Fokker-Planck equations with the novel perspective that solving the PDE can be reduced to independent instances of density estimation tasks based on the trajectories sampled from its associated particle dynamics. With this approach, one sidesteps error accumulation occurring from integrating the PDE dynamics on a parameterized function class. This approach significantly simplifies deployment, as one is free of the challenges of implementing loss terms based on the differential equation. In particular, we introduce a novel class of high-dimensional functions called the functional hierarchical tensor (FHT). The FHT ansatz leverages a hierarchical low-rank structure, offering the advantage of linearly scalable runtime and memory complexity relative to the dimension count. We introduce a sketching-based technique that performs density estimation over particles simulated from the particle dynamics associated with the equation, thereby obtaining a representation of the Fokker-Planck solution in terms of our ansatz. We apply the proposed approach successfully to three challenging time-dependent Ginzburg-Landau models with hundreds of variables.
\end{abstract}

\end{frontmatter}

%========================================
\section{Introduction}

Solving high-dimensional partial differential equations (PDE) remains one of the main challenges in scientific computing, with applications in molecular dynamics, quantum mechanics, and high-dimensional control theory. A key element underlying the challenges in solving high-dimensional PDE is the fact that the computational resources needed to carry out traditional numerical techniques often scale exponentially in the dimension of the problem \cite{han2018solving}. This paper focuses on solving the Fokker-Planck equation:
\begin{equation}\label{eqn: fokker-planck}
    \partial_{t}p = \nabla \cdot (p \nabla V) + \frac{1}{\beta}\Delta p, \quad x \in \Omega \subset \R^{d}, t \in [0, T],
\end{equation}
which governs the distribution function $p(t,x)$ of the particle system governed by the Langevin dynamics $dX_t = -\nabla V(X_t) dt + \sqrt{2\beta^{-1}}dB_{t}$ over a potential function $V$ with inverse temperature $\beta$.

Solving Fokker-Planck equations in low dimensions is a well-developed topic. This paper is concerned with the case of hundreds or even thousands of dimensions. One prevalent type of high-dimensional Fokker-Planck equation comes from the discretization of an infinite-dimensional functional equation into finite but very large dimensions.  Throughout this text, we use the motivating example of a time-dependent Ginzburg-Landau (G-L) model used for studying the phenomenological theory of superconductivity \cite{ginzburg2009theory, hoffmann2012ginzburg, hohenberg2015introduction,weinan2004minimum}.  In the 2D G-L model, one can intuitively think of a particle taking the form of a field \(x(a) \colon [0,1]^2 \to \R\), where the potential function is a functional $V(x)$ defined as follows:
\begin{equation}\label{eqn: GL potential infinite}
  V(x) = \frac{\lambda}{2}\int_{[0,1]^2} |\nabla_a x(a)|^2 \, da + \frac{1}{4\lambda}\int_{[0,1]^2} |1- x(a)^2|^2 \, da,
\end{equation}
where \(\lambda\) is a parameter balancing the relative importance of the first term and the second term in \eqref{eqn: GL potential infinite}.
In the numerical treatment, one considers the discretization of the unit square \([0,1]^2\) into a grid of \(d = m^2\) points $\{(ih,jh)\}$, for
$h=\frac{1}{m+1}$ and $1\le i,j\le m$. The discretization of the field at this grid is the $d$-dimensional vector $x=(x_{(i,j)})_{1\le i,j \le m}$, where \(x_{(i, j)} = x(ih, jh)\). Under the discretization scheme, the potential energy is defined as
\begin{equation}\label{eqn: 2D GZ model}
V(x) = V(x_{(1,1)}, \ldots ,x_{(m,m)}) := \frac{\lambda}{2} \sum_{v \sim w}\left(\frac{x_{v} - x_{w}}{h}\right)^2 +  \frac{1}{4\lambda} \sum_{v}\left(1 - x_v^2\right)^2,
\end{equation}
where $v$ and $w$ are Cartesian grid points and \(v \sim w\) if and only if they are adjacent. The evolution of the distribution over discretized fields satisfies the Fokker-Planck equation \eqref{eqn: fokker-planck} with the discretized potential function $V(x)$. Even a moderate grid resolution results in a model very large in dimension, which is a recurring theme for discretizing infinite-dimensional functional differential equations. The associated functional differential equation is the equation for the functional $p(t,x)$ encoding the probability density at the function $x$ in time $t$.

%One can show that the evolution equation for \(P\) is the Wasserstein gradient flow of the following form
%\[\partial_{t} P = \mathrm{div}_{u}\left(P \cdot \nabla_{u} \frac{\delta E}{\delta P}\right),\]
%where \(\frac{\delta E}{\delta P}\) is the Frechet derivative of \(E\) with respect to \(P\), and the energy function \(E\) is defined by \(E(P) = \int_{u}V(u)P(u) \, du + \beta^{-1}\int_{u}P(u) \ln P(u) \, du\).

% \paragraph{Challenges of traditional numerical schemes}
\subsection{Background and contribution}
\paragraph{Solving Fokker-Planck equations without time-stepping}

Traditionally, numerical solvers for the high-dimensional Fokker-Planck equations rely on performing time integration of a given initial condition \(p|_{t = 0}\) along a chosen parametric function class. The time integration can be done either by directly performing the time stepping of the solution or it can be done by including the PDE loss in the regression target. A common feature of the time integration scheme is error accumulation, as errors from previous time steps tend to build up.
Moreover, compressing the differential operator $\nabla \cdot (p \nabla V)$ in \eqref{eqn: fokker-planck} is often performed with heuristics, leading to quite significant implementation challenges both in the design of the compression scheme and in its development into code. 

The key insight in this paper is that solving the high-dimensional time-dependent Fokker-Planck equation can be reduced to a series of density estimation tasks. The reasoning is rather simple: one can approximate the Fokker-Planck solution \(p(t,x)\) at any time \(t \in [0, T]\) through density estimation over the sampled trajectory at time \(t\). Thus, one can approximate the solution \(p(t,x)\) at a collection of time steps \(0 = t_0 \leq \cdots \leq t_{K} = T\), and the continuous-in-time solution \(p(t,x)\) for an arbitrary \(t \in [0, T]\) can be approximated through appropriate interpolation. More specifically, we take advantage of the associated particle dynamics to the Fokker-Planck equation:
\begin{equation}\label{eqn: sde for fokker-planck}
    dX_t = -\nabla V(X_t) dt + \sqrt{2\beta^{-1}}dB_{t}.
\end{equation}
By simulating the dynamics of \eqref{eqn: sde for fokker-planck} through stochastic dynamic equation (SDE) simulation, we obtain a particle-based empirical approximation of \(p(t,x)\) by independently sampled trajectory data \(\{X_i(t)\}_{t \in [0, T]}\) for \(i = 1, \ldots, N\), which constitutes the input to the density estimation task. Thus, if $p(t,x)$ belongs to a function class with an efficient denoising or density estimation procedure, we can approximate the Fokker-Planck solution via denoising this empirical distribution. In effect, this approach sidesteps both the error accumulation and the coding challenges of a time integration scheme. 

Here, we choose the functional hierarchical tensor as the ansatz of the PDE solution $p(t,x)$. The rationale for this choice is that the diffusive nature of the Fokker-Planck equations tends to keep the tensor rank under control. Compared to most neural network algorithms, the proposed ansatz can efficiently calculate the normalization constant (see Section \ref{sec: related} for detailed discussions). Moreover, while directly modeling the density through nonparametric density estimation suffers from the curse of dimensionality \cite{hastie2009elements}, our proposed approach effectively performs by denoising the empirical distribution through restriction to the proposed parametric function class, thus avoiding the curse of dimensionality. Among the possible choices for functional tensor networks, while an alternative is to perform density estimation with a functional tensor train \cite{chen2023combining, ren2023high, dektor2021dynamic, dektor2021rank, soley2021functional}, we choose a functional hierarchical tensor to better account for low-rank structures of more general high-dimensional Fokker-Planck equations.

\paragraph{Functional hierarchical tensor}

\begin{figure}[t!]
    \centering
    \includegraphics[width = 0.7\textwidth]{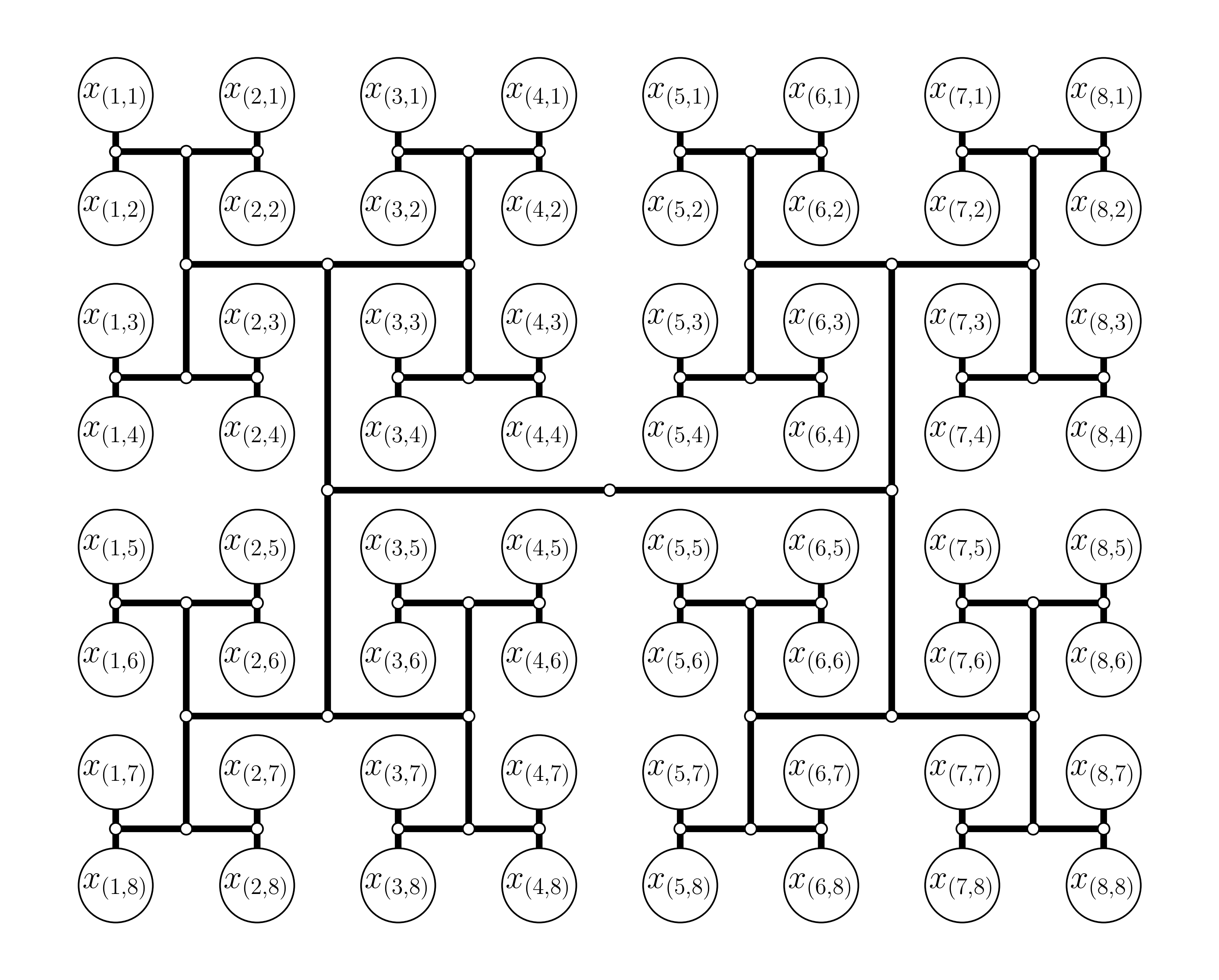}
    \caption{Illustration of the hierarchical bipartition of a 2D Ginzburg-Landau model discretized to an $8\times 8$ grid.}
    \label{fig:2D_GZ_bipar}
\end{figure}

Functional hierarchical tensor is understood in the literature  
\cite{bigoni2016spectral,gorodetsky2019continuous} to be a suitable tool for reduced-order modeling due to its versatility for multidimensional systems. In particular, the ansatz uses a modeling assumption natural for modeling discretized functional differential equations. The hierarchical tensor network assumes a low-rank structure along a hierarchical bipartition of the variables, which readily applies to discretized infinite-dimensional models. For example, 
the 2D G-L model has the natural geometry of a 2D grid, and it is natural to construct hierarchical bipartition of the variables through alternatively partitioning the variables along the two axes of the 2D grid, as shown in Figure \ref{fig:2D_GZ_bipar}. In comparison, it is known that the tensor train ansatz assumes the variables \(x = (x_1, \ldots, x_d)\) satisfy low-rankness according to the 1D structure set by the tensor network \cite{hur2023generative}, which makes it unsuitable for the 2D Ginzburg-Landau model, and the same argument applies to models for which the variables do not obey a 1D topological structure.

Our choice of functional hierarchical tensor over functional tensor train is partially motivated by parallel developments in quantum physics. For the task of modeling quantum systems, the tensor train ansatz is more commonly known as matrix product states \cite{ostlund1995thermodynamic}, where it is well-known that it is more suitable for modeling the ground state of a 1D system \cite{cirac2021matrix}. For modeling systems with nonlocal interaction, a hierarchical tensor network is often more advantageous \cite{vidal2007entanglement}. Through the use of a functional tensor network, our work extends the capability of the hierarchical tensor network into modeling continuous distributions.

\paragraph{Density estimation through hierarchical sketching}
After the trajectories are sampled through particle simulations, we process these trajectories through a sketching-based density estimation algorithm, obtaining a hierarchical tensor representation of the Fokker-Planck equation solution. The main contribution of this work is an end-to-end procedure that produces a Fokker-Planck solution from sampled trajectories. 

To ensure that an accurate approximation can be obtained with a small parameter size, we have worked out practical choices over key design questions. The important ingredients are the structure of the hierarchical bipartition, the sketch function, and the basis to be used for the hierarchical tensor. Using sophisticated numerical techniques, we are able to obtain solutions to Ginzburg-Landau models of very high dimensions, including solving the multidimensional Ginzburg-Landau model in \eqref{eqn: 2D GZ model} with \(d\) in the hundreds.

\subsection{Related work}\label{sec: related}

\paragraph{Neural network for solving high-dimensional PDEs} Neural network has shown tremendous success in recent years in solving high-dimensional PDEs. The most common approach to solve PDEs with neural networks primarily relies on penalizing the partial-differential equation either in terms of the strong form (as in the physics-informed neural network \cite{raissi2019physics} and the deep Galerkin method \cite{sirignano2018dgm}) or the variational form (as in the deep Ritz method \cite{yu2018deep}). Another separate direction is to solve the PDE by creating regression targets special to the problem itself, as can be seen in the deep BSDE method \cite{han2018solving} and neural operators \cite{kovachki2021neural}. While these algorithms perform well, further improving their effectiveness in high-dimensional PDEs is challenging.

\paragraph{Neural network methods for density estimation}
As our work primarily focuses on using density estimation to solve the Fokker-Planck equation, our work is related to generative learning algorithms for which density estimation is possible. Examples include normalizing flows \cite{tabak2010density,rezende2015variational}, energy-based models \cite{hinton2002training, lecun2006tutorial}, and the more recent diffusion models \cite{song2019generative,song2021maximum}. Only normalizing flows allow direct density evaluation. In contrast, the energy-based model and diffusion model only allow for approximating the likelihood function, which renders them inapplicable in our case, as the objective is to obtain the model's associated normalized likelihood function. Moreover, training several neural networks across many time steps might be inefficient, while our sketching-based approach relies on quick linear algebraic subroutines and does not involve iterative optimization.

\paragraph{Monte-Carlo method} 
As one has access to trajectories of the particle dynamics due to Monte-Carlo sampling, one might hope to perform density estimation through nonparametric density estimation directly. However, particle-based density estimations such as kernel density estimation suffer from the curse of dimensionality \cite{hastie2009elements}, making it ill-suited to our problem. If the goal is to estimate the statistical moments of the solved Fokker-Planck equation, then indeed the traditional Monte-Carlo approach is well-suited to estimate such quantity with an \(O(1/\sqrt{N})\) rate \cite{liu2001monte}, \(N\) being the number of samples. We remark that in practical cases, the variance is often too high for the estimation to have a satisfactory estimation accuracy. 

\paragraph{Tensor network for high-dimensional PDE}
There are several works using tensor networks to solve high-dimensional PDEs. One approach is to use grid discretization of the PDE solution \(p\), resulting in the task of solving the PDE for a \(d\)-dimensional tensor, as can be seen in \cite{chen2023committor}. The functional tensor network approach employs a mesh-free representation of the PDE solution based on a tensor network. Along this line of work, the considered tensor network structures include canonical-polyadic (CP) \cite{khoromskij2011tensor} and tensor train \cite{bigoni2016spectral,gorodetsky2019continuous,soley2021functional,dektor2021dynamic,eigel2019non}. The use of hierarchical tensor in solving high-dimensional PDE has also been discussed \cite{bachmayr2016tensor}, though quite different from our proposed methodology. The work in \cite{richter2021solving,chen2023combining} bears the most resemblance to our work, as they combine particle methods with a functional tensor network structure. The main difference is the fact that our proposed algorithm does not involve a time-stepping component, and therefore, our method is less prone to error accumulation. A similar sketching-based density estimation algorithm for solving the Fokker-Planck equation with a functional tensor train can be done by combining our approach with the sketching methods outlined in \cite{hur2023generative}.

\subsection{Contents and notations} \label{sec: notation}

We outline the structure of the remainder of the manuscript. Section \ref{sec: background} gives a detailed introduction to functional hierarchical tensor and its sketching-based density estimation algorithm. Section \ref{sec: alg} goes through the algorithm for solving the Fokker-Planck equation with the functional hierarchical tensor. Section \ref{sec: numerics} details the numerical implementations and results from solving Ginzburg-Landau models in multiple dimensions.

For notational compactness, we introduce several shorthand notations for simple derivation. For \(n \in \mathbb{N}\), let \([n] := \{1,\ldots, n\}\). For an index set \(S \subset [d]\), we let \(x_{S}\) stand for the subvector with entries from index set $S$. We use \(\Bar{S}\) to denote the set-theoretic complement of $S$, i.e. \(\Bar{S} = [d] - S\).

%========================================
\section{Hierarchical functional tensor network}\label{sec: background}

To address the Ginzburg-Landau and general models coming from discretizing high-dimensional functional equations, we propose to use the function class of a functional hierarchical tensor based on a binary-tree-based low-rank structure.
In this section, we go over the functional hierarchical tensor and its sketching algorithm. The symbol \(d\) is reserved for the dimension of the state-space variable, and the symbol \(N\) is reserved for the number of samples. Without loss of generality, we shall go over the procedure to perform density estimation over a collection of \(N\) samples \(\{y^{(i)} \in \R^d\}_{i = 1}^{N}\). One can think of \(y^{(i)}\) as the data collected from the \(i\)-th trajectory at a fixed time \(t \in [0, T]\). 

\subsection{Hierarchical bipartition}

% Describe what this is. Tree structure, ordering. Basis function at the leaf.

We first describe the functional tensor network representation of a \(d\)-dimensional function \(p \colon \mathbb{R}^{d} \to \mathbb{R}\). 
In general, let \(\{\psi_{i}\}_{i = 1}^{n}\) denote a collection of orthonormal function basis over a single variable, and let \(C \in \R^{n^d}\) be the tensor represented by a tensor network. The \emph{functional tensor network} is the \(d\)-dimensional function defined by the following equation:
\begin{equation}\label{eqn: htn forward map}
    p(x)\equiv p(x_{1}, \ldots, x_{d}) = \sum_{i_{1}, \ldots, i_{d} = 0}^{n-1} C_{i_1,\ldots, i_d} \psi_{i_1}(x_1)\cdots \psi_{i_d}(x_d) = \left<C, \, \bigotimes_{j=1}^{d} \Vec{\Psi}(x_j) \right>,
\end{equation}
where \(\Vec{\Psi}(x_j) = \left[\psi_{1}(x_j),\ldots, \psi_{n}(x_j) \right]\) is an \(n\)-vector encoding the evaluation of any \(x_j\) over the entire single variable function basis set. 

In this work, we choose \(C\) to be represented by a hierarchical tensor network ansatz so that the associated \(p\) takes the form of a functional hierarchical tensor. In addition to density evaluation, this ansatz allows for efficient evaluation of moments and efficient sampling (see Section \ref{sec: application} for more details on potential application).

\begin{figure}[h!]
    \centering
    \begin{subfigure}{0.45\textwidth}
        \centering
        \includegraphics[width=\textwidth]{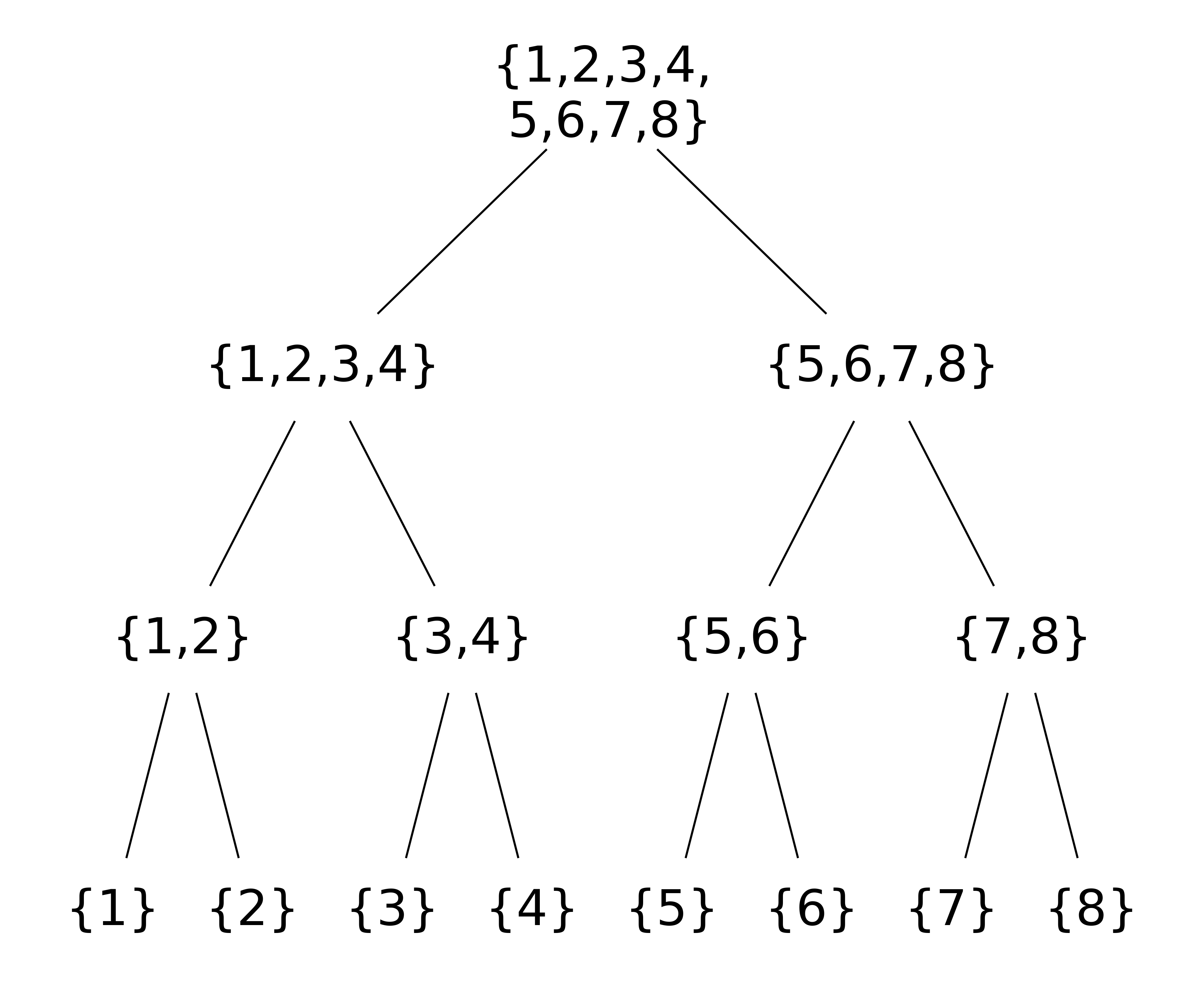} % first figure itself
        \caption{Binary decomposition of variables}
        \label{fig:binary_tree_8_nodes_subfig}
    \end{subfigure} \hfill % <-- added \hfill
    \begin{subfigure}{0.54\textwidth}
        \centering
        \includegraphics[width=\textwidth]{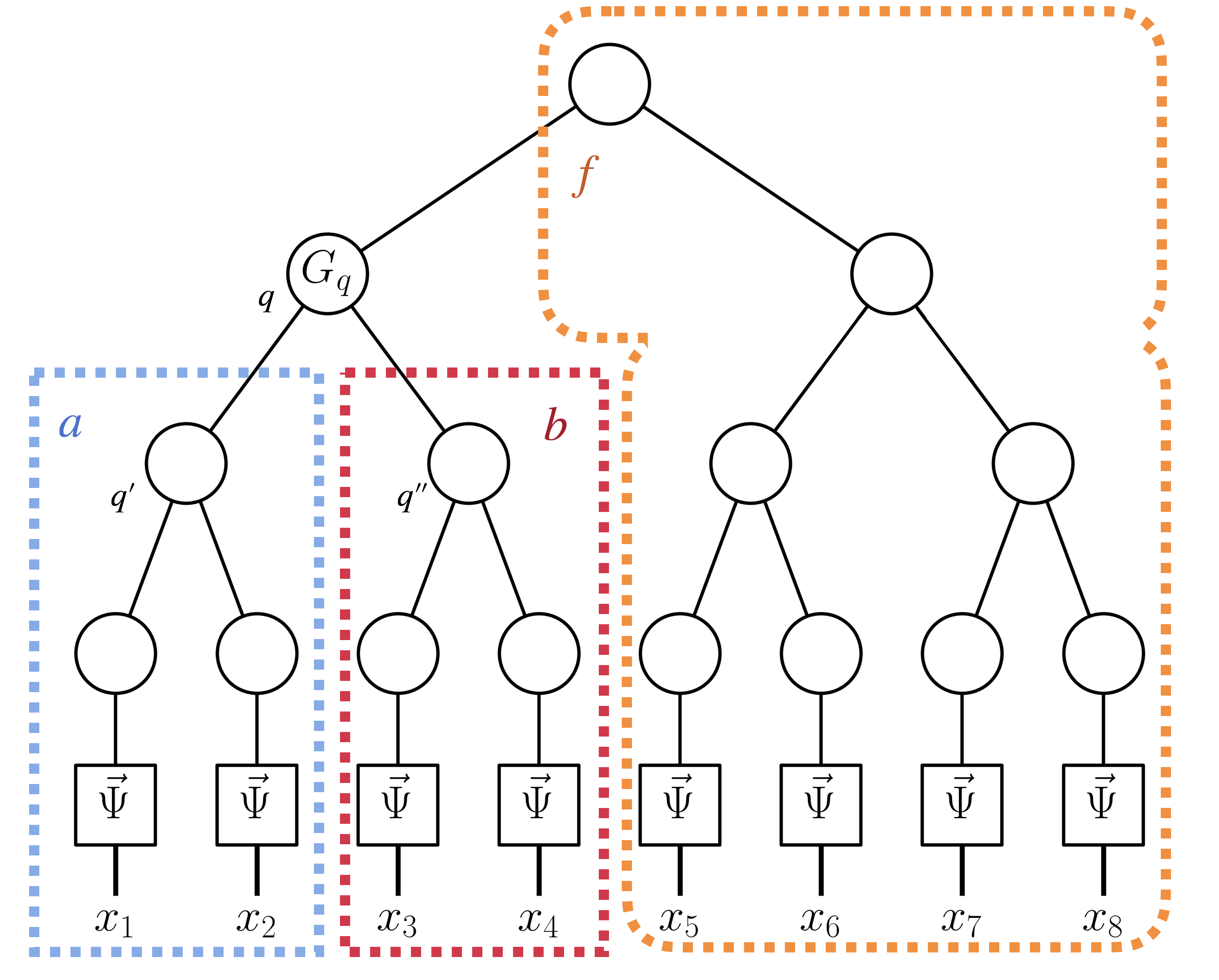} % second figure itself
        \caption{Functional hierarchical tensor diagram}
        \label{fig:FHT_L_3}
    \end{subfigure}
    \caption{Illustrations of functional hierarchical tensor with $d = 8$.}
    \label{fig:decompositions}
\end{figure}

Crucial to the notion of the hierarchical tensor network is a hierarchical bipartition of the variable set. Without loss of generality, let \(d = 2^{L}\) so that the variable set admits exactly \(L\) levels of variable bipartition. As illustrated in Figure \ref{fig:binary_tree_8_nodes_subfig}, at the \(l\)-th level, the variable index set is partitioned according to 
\begin{equation}\label{eqn: bipartition}
    [d] = \bigcup_{k = 1}^{2^{l}} I_{k}^{(l)}, \quad I_{k}^{(l)} := \{ 2^{L - l}(k-1) + 1, \ldots, 2^{L - l}k\},
\end{equation}
which in particular implies the recursive relation that \(I_{k}^{(l)} = I_{2k-1}^{(l+1)}\cup I_{2k}^{(l+1)}\). 

\subsection{Tensor network structure}
We introduce the fundamental structure of a functional hierarchical tensor. 
Importantly, we focus on a single node \(q\), say at level index \(l\) and block index \(k \in [2^{l}]\). We take the shorthand \(a(q):= I_{2k-1}^{(l+1)}, b(q):= I_{2k}^{(l+1)}\) and \(f(q) = [d] - a(q) \cup b(q)\), to be its left branch, right branch, and the rest, respectively, as illustrated in Figure \ref{fig:FHT_L_3}. While these index sets depend on the node $q$, in order to simplify the formulae, we refer to them just as $a$, $b$, and $f$. 

Intrinsically, a hierarchical tensor network represents a tensor \(C \in \mathbb{R}^{n^d}\) with a low-rank structure along every hierarchical bipartition. The bipartition at the node \(q\) is defined to  be \([d] = \Bar{f} \cup f\), where \(\Bar{f}\) is the complement of \(f\) as a subset of \([d]\). With the partition \([d] = \Bar{f} \cup f\), there exists \(r_{f} \in \mathbb{N}\) so that the \(n^{|\Bar{f}|} \times n^{|f|}\) unfolding matrix \(C(i_{\Bar{f}}; i_f)\) is of rank \(r_{f}\). In other words, there exist \(C_{\Bar{f}} \colon \R^{n^{|\Bar{f}|}} \times \R^{r_{f}}\) and \(C_f \colon \R^{r_f} \times  \R^{n^{|f|}} \) that form a low-rank decomposition of \(C\) (illustrated in Figure \ref{fig:tensor_diagram_equation}(a)):
\begin{equation}\label{eqn: unfolding matrix of C}
    C(i_{1}, \ldots, i_{d}) = \sum_{\theta = 1}^{r_f}C_{\Bar{f}}(i_{\Bar{f}}, \theta)C_f(\theta, i_{f}).
\end{equation}

If the low-rankness property in \eqref{eqn: unfolding matrix of C} holds for any node \(q\), the exponential-sized tensor \(C\) can be described by a mere \(O(dr^3)\) number of parameters by a hierarchical tensor network, where \(r:= \max_{q}(r_{f(q)})\). The tensor network is parameterized by a collection of tensor cores denoted by \(\{G_{q}\}_{q}\).
As Figure \ref{fig:decompositions} illustrates, the tensor cores have the same tree structure as in the variable bipartition. Each internal bond of 
the tensor network is associated with a bond in the binary decomposition of the variables. Each physical bond of the tensor network connects a leaf node with a $\Vec{\Psi}(x_j)$ node. In summary, the graphical structure of the hierarchical tensor network is formed by a binary tree, and the physical index is at the leaf node. The fact that such cores exist can be proved by induction using \eqref{eqn: unfolding matrix of C}, and a proof can be found in \cite{peng2023generative}. In terms of runtime complexity and memory complexity, the functional hierarchical tensor enjoys the same \(O(d)\) scaling as that of a hierarchical tensor \cite{peng2023generative}.

\subsection{Sketching algorithm}\label{sec: sketch}

% Under this choice, one can decompose the

% Describe the hierarchical sketching algorithm here.

% While the sketching algorithm for hierarchical tensor for discrete distributions has been developed \cite{peng2023generative}, one of our technical novelty is the extension into functional tensor networks \cite{bigoni2016spectral, oseledets2013constructive, gorodetsky2019continuous}.

\begin{figure}[h!]
    \centering
    \includegraphics[width = 0.8\textwidth]{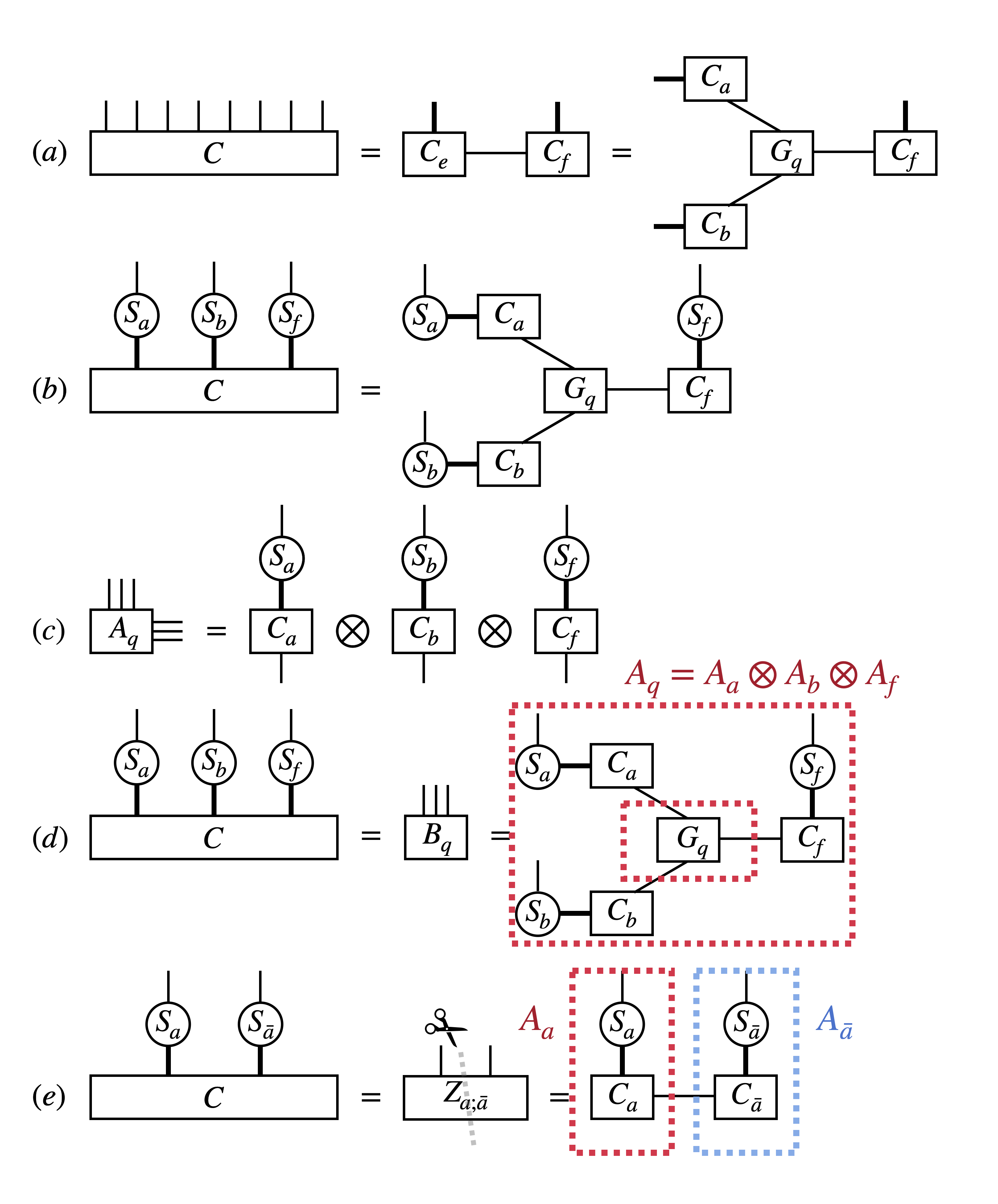}
    \caption{Tensor diagram representation of main equations used in Section \ref{sec: sketch}. Equation \eqref{eqn: unfolding matrix of C} and \eqref{eqn: unsketched linear system for G} are summarized in (a). Equation \eqref{eqn: sketched linear system for G} is shown in (b). The coefficient term in \eqref{eqn: sketched linear system for G} is shown in (c). The definition of \(B_q\) is shown in (d). Equation \eqref{eqn: eq for Z_ag} for \(Z_{a;\Bar{a}}\) is shown in (e), where the scissor symbol indicates the rank-\(r_{a}\) compression through SVD. }
    \label{fig:tensor_diagram_equation}
\end{figure}

We shall go over how to use a sketch-based method to obtain a functional hierarchical tensor representation from the given collection of samples \( \left\{ y^{(i)}:=\left(y_1^{(i)}, \ldots, y_{d}^{(i)} \right) \right\}_{i = 1}^{N}\). For the reader's convenience, important equations in this subsection are included in Figure \ref{fig:tensor_diagram_equation} in terms of the tensor diagram. Below, we review the main equation behind hierarchical sketching in the functional case. The goal is to solve for the tensor core \(G_{q}\), where \(q\) is the node of the \(k\)-th block at level \(l\).

% Under the assumption that the density function \(p^{\star}\) satisfies the functional hierarchical tensor network structure in \eqref{eqn: htn forward map}, one can use randomized linear algebra techniques to estimate each tensor component consistently. We remark that the sketching algorithm for hierarchical tensor for discrete distributions has been developed \cite{peng2023generative}. Conceptually, our proposed algorithm extends the hierarchical tensor sketching in \cite{peng2023generative} into the functional case. We remark that the connection is clear in the following sense: By calculating the statistical moments of the samples systematically, one can obtain a sketched linear system for each tensor component in the hierarchical tensor \(C\), which is agnostic of whether \(C\) encodes discrete or continuous densities. 

\paragraph{Equations} For simplicity, we first assume $0< l < L$ so that \(q\) is neither the root nor the leaf node. Similar as before, we have \(a := I_{2k-1}^{(l+1)}, b := I_{2k}^{(l+1)}\) and \(f = [d] - a \cup b\). Then, the structural low-rankness property of \eqref{eqn: unfolding matrix of C} implies that there exist \(r_a, r_b, r_f \in \mathbb{N}\), \(C_a \colon [n^{|a|}] \times [r_a] \to \R\), \(C_b \colon [n^{|b|}] \times [r_b] \to \R\), \(C_f \colon [r_f] \times [n^{|f|}] \to \R\) such that the following equation holds (illustrated in Figure \ref{fig:tensor_diagram_equation}(a)):
\begin{equation}\label{eqn: unsketched linear system for G}
    C(i_{1}, \ldots, i_d) = \sum_{\alpha, \beta, \theta}C_{a}(i_{a}, \alpha)C_{b}(i_{b}, \beta)G_{q}(\alpha, \beta, \theta)C_{f}(\theta, i_{f}).
\end{equation}

Equation \eqref{eqn: unsketched linear system for G} is exponential-sized, and therefore, one would not solve this equation directly.
The hierarchical sketching algorithm essentially solves the over-determined linear system \eqref{eqn: unsketched linear system for G} for \(G_{q}\) with the use of sketch functions. Let \(\tilde{r}_a, \tilde{r}_b, \tilde{r}_f\) be integers so that \(\tilde{r}_a > r_a, \tilde{r}_b > r_b, \tilde{r}_f > r_f\). Through the sketch functions \(S_a \colon [n^{|a|}] \times [\tilde{r}_a] \to \R, S_b \colon [n^{|b|}] \times [\tilde{r}_b] \to \R, S_f \colon [n^{|f|}] \times [\tilde{r}_f] \to \R\), we can contract the linear system in \eqref{eqn: unsketched linear system for G} with \(S_a, S_b, S_f\) at the variables \(i_a, i_b, i_f\) respectively, which leads to the following sketched linear system (illustrated in Figure \ref{fig:tensor_diagram_equation}(b) and Figure \ref{fig:tensor_diagram_equation}(d)):
\begin{equation}\label{eqn: sketched linear system for G}
    B_q(\mu, \nu, \zeta) = \sum_{\alpha, \beta, \theta}A_{a}(\mu, \alpha)A_b(\nu, \beta)A_f(\zeta, \theta)G_{q}(\alpha, \beta, \theta),
\end{equation}
where \(A_{a}, A_b, A_f\) are respectively the contraction of \(C_a, C_b, C_f\) by \(S_a, S_b, S_f\) and \(B_q\) is the contraction of \(C\) by \(S_a \otimes S_b \otimes S_f\), as illustrated by the first equality of Figure \ref{fig:tensor_diagram_equation}(d). The equation \eqref{eqn: sketched linear system for G} can be seen as the following linear system of \(\tilde{r}_a\tilde{r}_b\tilde{r}_f\) equations for the unknown \(G_{q}\) (illustrated in Figure \ref{fig:tensor_diagram_equation}(d)):
\[
(A_{a} \otimes A_{b} \otimes A_{f})G_{q} = B_{q}.
\]

% Figure \ref{fig:tensor_diagram_equation} gives the tensor diagram for \eqref{eqn: sketched linear system for G} as well as the tensor diagram for the tensor \(B_q\).

While obtaining \(B_q\) in \eqref{eqn: sketched linear system for G} is done by contracting known tensors, obtaining terms such as \(A_a\) is only possible after specifying the gauge degree of freedom for \(C_a\). However, we show that one can obtain \(A_a\) without having access to \(C_a\). The structural low-rankness of the ansatz implies the existence of \(C_{\Bar{a}}\) such that there exists a linear system \(C(i_1, \ldots, i_d) = \sum_{\alpha}C_a(i_a, \alpha)C_{\Bar{a}}(\alpha, i_{{\Bar{a}}})\).
One can sketch this linear system by contracting it with \(S_a, S_{{\Bar{a}}}\), which leads to the following linear system (illustrated in Figure \ref{fig:tensor_diagram_equation}(e)):
\begin{equation}\label{eqn: eq for Z_ag}
    Z_{a;{\Bar{a}}}(\mu, \psi) = \sum_{\alpha}A_a(\mu, \alpha)A_{{\Bar{a}}}(\alpha, \psi),
\end{equation}
where \(A_{{\Bar{a}}}\) is the contraction of \(C_{{\Bar{a}}}\) with \(S_{{\Bar{a}}}\). Then, due to a gauge degree of freedom in choosing \(C_a\), there also exists a gauge degree of freedom in choosing \(A_a\). Thus, one can perform singular value decomposition (SVD) on \(Z_{a;{\Bar{a}}}\), obtaining \(Z_{a;{\Bar{a}}}(\mu, \psi) = \sum_{\alpha}U(\mu, \alpha)V(\alpha, \psi)\), and the choice of \(A_a = U\) and \(A_{{\Bar{a}}}= V\) forms a consistent choice of gauge between the pair \((C_a, C_{{\Bar{a}}})\). 

Likewise, one can obtain \(Z_{b; \Bar{b}}\) and the SVD of \(Z_{b;\Bar{b}}\) results in \(A_b\) and \(A_{\Bar{b}}\). The same holds for
 \(Z_{\Bar{f};f}\), \(A_f\), and \(A_{\Bar{f}}\). To solve for \(G_{q}\) in the linear system \eqref{eqn: sketched linear system for G}, one simply contracts \(B_q\) with the pseudo-inverse of \(A_a, A_b, A_f\). Thus, one can use \eqref{eqn: sketched linear system for G} to solve for \(G_{q}\) for any \(l \not = 0, L\).

Our construction likewise gives rise to the equation of \(G_{q}\) for \(l = 0, L\) as special cases. For \(l = 0\), one can go through with the same calculation by setting \(C_f= S_f = A_f = 1\) in \eqref{eqn: unsketched linear system for G} and \eqref{eqn: sketched linear system for G}.  For \(l = L\), one can sketch the linear system in \eqref{eqn: unsketched linear system for G} for \(f = [d] - \{k\}\) by contraction with \(S_{f}\). The detail for both cases is described in Algorithm \ref{alg:1}.

\paragraph{Approximation via samples}
We now describe how one can obtain a consistent estimation of \(G_{q}\) through samples. As before, suppose \(\{\psi_{i}\}_{i = 1}^{n}\) is the collection of univariate orthonormal function basis which defines the functional tensor network basis. To the sketch tensor \(S_a \colon [n^{|a|}] \times [\tilde{r}_a] \to \R \), one associates a continuous sketch function \(s_a \colon \R^{|a|} \times [\tilde{r}_a] \to \R\) through the following construction:
\begin{equation}\label{eqn: sketch function}
    s_a(x_{a}, \mu) = \sum_{i_{j}, j \in a} S_a(i_{a}, \mu) \prod_{j \in a}\psi_{i_j}(x_j),
\end{equation}
and one can likewise obtain \(s_{\Bar{a}}\) through \(S_{\Bar{a}}\).
By orthonormality of the function basis, it is easy to check that
\[
Z_{a; \Bar{a}}(\mu, \psi) = \int_{\R^d}  p(x_1, \ldots, x_d) s_a(x_{a}, \mu) s_{\Bar{a}}(x_{\Bar{a}}, \psi)\, d \, x_1 \ldots x_d = \mathbb{E}_{X \sim p}\left[s_a(X_{a}, \mu) s_{\Bar{a}}(X_{\Bar{a}}, \psi)\right],
\]
and the formulae for $(Z_{b; \Bar{b}}$ and $Z_{\Bar{f}; f}$ follow likewise. Similarly, one can obtain \(B_q\) in \eqref{eqn: sketched linear system for G} by 
\[
B_q(\mu, \nu, \zeta) = \mathbb{E}_{X \sim p}\left[s_a(X_{a}, \mu) s_b(X_b, \nu) s_{f}(X_f, \zeta)\right].
\]

Since \(\{y^{(i)} \in \R^d\}_{i = 1}^{N}\) provides an empirical approximation to $p$, one obtains a finite sample approximation of \(Z_{a; \Bar{a}}\) in the following way:
\begin{equation}\label{eqn: Z_ag approximation}
    Z_{a; \Bar{a}}(\mu, \psi) \approx \frac{1}{N} \sum_{i = 1}^{N} s_a(y^{(i)}_{a}, \mu) s_{\Bar{a}}(y^{(i)}_{\Bar{a}}, \psi)
\end{equation}
and likewise one can also obtain a sample estimation of \(B_q\) and thus approximately solve for \(G_{q}\). It is easy to check that a consistent estimator of \(Z_{a;\Bar{a}}\) leads to a consistent estimator of \(A_a\), and therefore one can likewise consistently estimate \(A_b, A_f\), and thus the approximated solution is a consistent estimator for \(G_{q}\). 

\paragraph{Algorithm summary}
The sketching algorithm is summarized as Algorithm \ref{alg:1}, which demonstrates the whole procedure to carry out the functional hierarchical tensor sketching algorithm given a sample collection, including the details for the omitted edge cases of \(l = 0, L\).

\begin{algorithm}[h]
\caption{Functional hierarchical tensor sketching.}
\label{alg:1}
\begin{algorithmic}[1]
\REQUIRE Sample \(\{y^{(i)}\}_{i = 1}^{N}\).
\REQUIRE Chosen function basis \(\{\psi_i\}_{i = 1}^{n}\).
\REQUIRE Collection of sketch functions \(\{s_{I_{k}^{(l)}}\}, \{s_{[d] - I_{k}^{(l)}}\}\) and target internal ranks $\{r_{I_{k}^{(l)}}\}$ for each level index \(l\) and block index \(k\).
\FOR{each node \(q\) on the hierarchical tree}
    \STATE Set \(l\) as the level index of \(q\). Set \(k\) as the block index of \(q\).
    \IF{\(q\) is not leaf node}
        \STATE \((a, b) \gets (I_{2k-1}^{(l+1)}, I_{2k}^{(l+1)})\)
        \STATE Obtain \(Z_{a; \Bar{a}}, Z_{b; \Bar{b}}\) by \eqref{eqn: Z_ag approximation}
        \STATE Obtain \(A_{a}\) as the left factor of the best rank \(r_{a}\) factorization of \(Z_{a; \Bar{a}}\)
        \STATE Obtain \(A_{b}\) as the left factor of the best rank \(r_{b}\) factorization of \(Z_{b; \Bar{b}}\). 
    \ENDIF
    \IF{\(q\) is not root node}
        \STATE \(f \gets [d] -I_{k}^{(l)}\)
        \STATE Obtain \(Z_{\Bar{f}; f}\) by \eqref{eqn: Z_ag approximation}. 
        \STATE Obtain \(A_{f}\) as the right factor of the best rank \(r_{f}\) factorization of \(Z_{\Bar{f}; f}\).
    \ENDIF
    \IF{\(q\) is leaf node}
    \STATE Obtain \(
        B_{q}(\mu, \nu) = \frac{1}{N}\sum_{i}s_a(y^{(i)}_a, \mu) s_b(y^{(i)}_b, \nu).
        \)
    \STATE Obtain \(G_{q}\) by solving the over-determined linear system \((A_{a} \otimes A_{b})G_{q} = B_{q}\).
    \ELSIF{\(q\) is root node}
    \STATE Obtain \(
        B_{q}(j, \zeta) = \frac{1}{N}\sum_{i}\psi_{j}(y^{(i)}_k) s_{f}(y^{(i)}_f, \zeta).\)
    \STATE Obtain \(G_{q}\) by solving the over-determined linear system \((A_{f})G_{q} = B_{q}\).
    \ELSE
    \STATE Obtain \(
        B_{q}(\mu, \nu, \zeta) = \frac{1}{N}\sum_{i}s_a(y^{(i)}_a, \mu) s_b(y^{(i)}_b, \nu) s_{f}(y^{(i)}_f, \zeta).
        \)
    \STATE Obtain \(G_{q}\) by solving the over-determined linear system \((A_{a} \otimes A_{b} \otimes A_{f})G_{q} = B_{q}\).
    \ENDIF
\ENDFOR
\end{algorithmic}
\end{algorithm}

%========================================
\section{Solving for Fokker-Planck equation}\label{sec: alg}\label{sec: application}

In this section, we detail our approach to solving the Fokker-Planck equation. The given information regarding the PDE includes the potential function \(V \colon \mathbb{R}^{d} \to \mathbb{R}\), terminal time \(T\), and initial distribution \(p_{0}\). The main workflow is summarized in the following simple steps:
\begin{enumerate}
    \item Select time steps \(0 = t_0 \leq t_1 \leq \ldots \leq t_{K} = T\). Then, with a chosen numerical scheme (e.,g. Euler-Maruyama), run \(N\) independent instance of SDE simulations on the stochastic dynamic equation \eqref{eqn: sde for fokker-planck} to time \(T\) with the initial condition \(p_{0}\). Obtain \(N\) trajectories and record the trajectory data \(\{X^{(i)}(t_{j})\}_{j \in [K], i \in [N]}\).

    \item For \(j = 1, \ldots, K\), obtain an approximated functional hierarchical solution \(p_j(x) \approx p(x,t_j)\) by calling Algorithm \ref{alg:1} on the sample \(\{X^{(i)}(t_{j})\}_{i \in [N]}\).
    
    \item Output the continuous-in-time solution for \(p(t,x)\), defined as follows: Suppose that \(t \in [0, T]\) satisfies \(t \in (t_j, t_{j+1})\), then \(p(t, x)\) can be approximated with \(p_{j}(x)\) and \(p_{j+1}(x)\) via an appropriate interpolation. 
\end{enumerate}

Once the approximated solution \(p\) is obtained, one can perform sampling or compute observables based on the ansatz. As \(p(t, x)\) interpolates between \(p_j\) and \(p_{j+1}\) for \(t \in (t_j, t_{j+1})\), the task of sampling and computing observables reduces to performing such tasks on snapshot solutions \(p_{j}\) and \(p_{j+1}\). Therefore, it suffices to illustrate how to perform such tasks on an obtained snapshot solution \(p_j\).

\begin{figure}[h!]
    \centering
    \includegraphics[width = 0.4\textwidth]{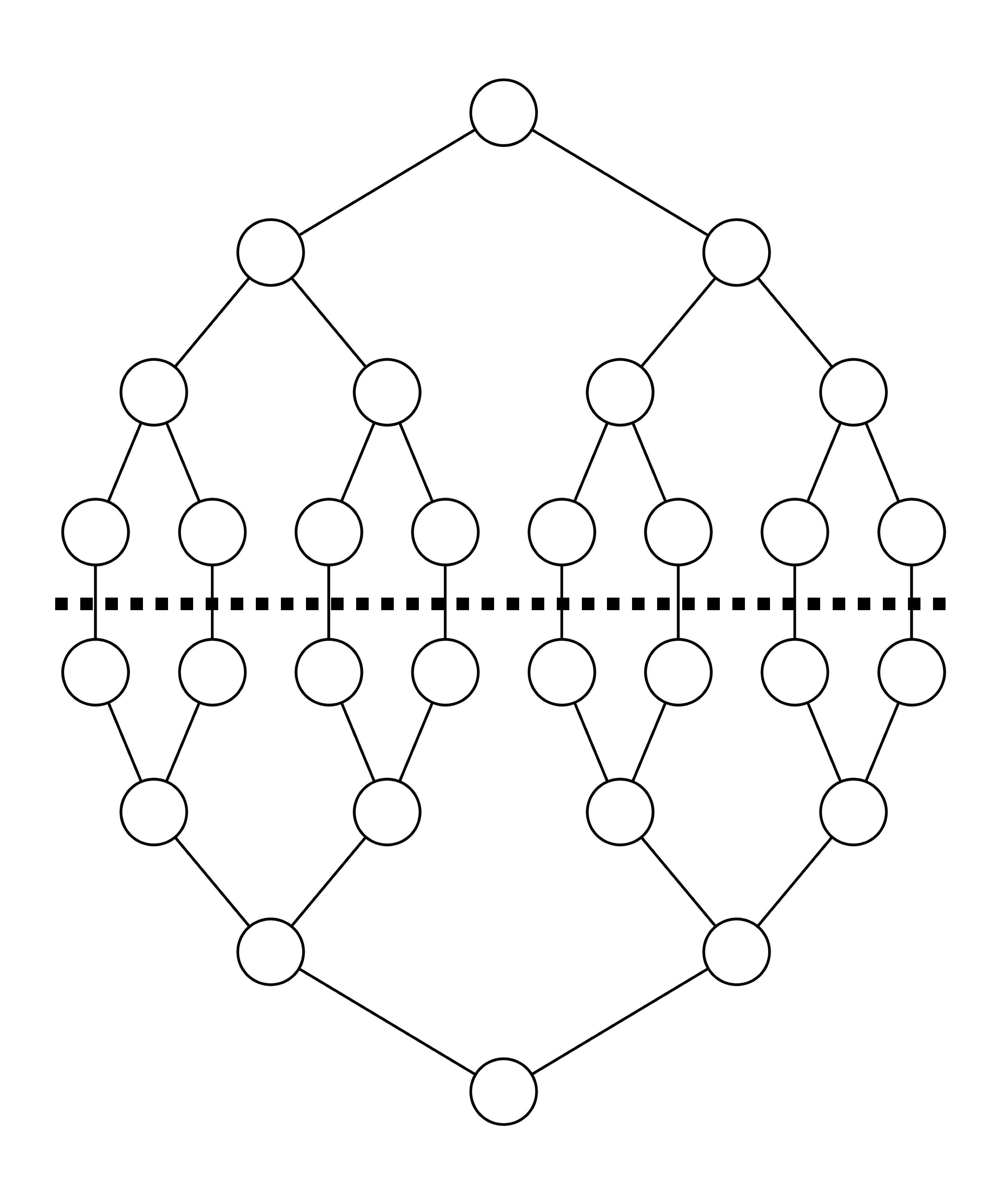}
    \caption{Tensor diagram representation for the measuring observable \(\mathbb{E}_{X \sim p_j}[O(X)]\). The top half of the tensor cores represents the tensor network for \(p_j\), and the bottom half represents the tensor network for \(O\).
    % The tensor cores for \(p_{j}\) is denoted by \(\{G_{k}^{(l)}\}_{k, l}\), and the tensor cores for \(O\) is denoted by \(\{F_{k}^{(l)}\}_{k, l}\).
    }
    \label{fig:Observable_obtain}
\end{figure}

\paragraph{Application I: Observable estimation}
For an observable function \(O \colon \R^d \to \R\), note that one has \(\mathbb{E}_{X \sim p_j}[O(X)] = \left<O(x), p_j(x)\right>_{L^2(\R^d)}\).
When \(O\) can be written in terms of a finite sum of a monomial over the basis functions, one can compute the inner product efficiently through a tensor contraction, as the orthonormality of the function basis implies the function space inner product will coincide with the inner product between the coefficient tensor.  In general, when \(O\) is a functional hierarchical tensor, one can compute \(\left<O(x), p_j(x)\right>_{L^2(\R^d)}\) through tensor diagram contraction, as is shown in Figure \ref{fig:Observable_obtain}. The tensor contraction diagram in Figure \ref{fig:Observable_obtain} admits efficient evaluation through iterative tensor contractions starting from the middle two layers.

%For \(U, V \subset [d]\), we use \(p(x_{U})\) to denote marginal distribution of \(p\) when restricted to \(x_{U}\), and we use \(p(x_{U} | x_{V})\) to denote conditional distribution of \(x_{U}\) given \(x_{V}\). Then,

\paragraph{Application II: Sample generation}

The counting-to-sampling \cite{jerrum2003counting} perspective reduces the sampling task to the previous one for observable estimation. Using conditional distributions, we have the following decomposition of \(p_{j}\) at time $t_j$
\[
p_j(x_1, \ldots, x_d) = p_j(x_1) \cdot p_j(x_2 | x_1) \ldots p_j(x_d | x_{1}, \ldots, x_{d-1}).
\]
One then uses sequential Monte Carlo (see \cite{liu2001monte} for details) to generate i.i.d. samples: For \(k \in [d]\), suppose one has sampled the distribution up to index \(k-1\) with \(x_{i} = a_{i}\) for \(i = 1, \ldots, k - 1\). One can sample the \(k\)-th variable through \(X_{k} \sim p_{j}(x_k = x | x_{[k-1]} = a_{[k-1]})\), where \(p_{j}(x_{k} = x | x_{k-1} = a_{k-1}, \ldots, x_{1} = a_1) = \mathbb{E}_{X \sim p_j}[\prod_{i \in [k - 1]}\delta(X_i - a_i)\delta(X_k - x))]\). 

\subsection{Extensions}

Though we have been focused on Langevin dynamics for simplicity, our method can be generalized to a wider range of settings. We shall remark a few notable extensions to the current approach. 

\paragraph{Fokker-Planck beyond Langevin dynamics}
First, one can readily see that the approach is agnostic to the simulated SDE structure considered. In particular, one can consider the generalized Fokker-Planck equation,
\[
\partial_t p = \nabla \cdot (p f 
) + \sum_{i, j, k =1}^{d} \partial_{x_i x_{j}} (\sigma_{ik}\sigma_{jk} p), \quad x \in \Omega \subset \R^{d}, t \in [0, T],
\]
for a general particle system with the following form:
\begin{equation}\label{eqn: general sde for fokker-planck}
    dX_t = -f(X_t, t) dt + \sigma(X_t, t)dB_{t},
\end{equation}
where \(f \colon \R^d \times [0, T] \to \R^d\) and \(\sigma \colon \R^d \times [0, T] \to \R^{d\times d}\) are general advection and diffusion terms. While Algorithm \ref{alg:1} does not depend on the exact method used to generate the sample trajectory, it might be advisable to use a higher-order SDE simulation to guarantee sample quality in simulating \eqref{eqn: general sde for fokker-planck}.

\paragraph{Non-linear Fokker-Planck equations}
Secondly, one can also extend our method to nonlinear Fokker-Planck equations. Using the perspective that the Fokker-Planck equation is the Wasserstein gradient flow on an Energy potential \cite{jordan1998variational}, one can consider a functional  
\[
E(p) + \beta^{-1}\int p(x) \ln p(x) dx
\]
with a non-linear $E(p)$. The associated PDE is
\[
\partial_t p(t,x) = \nabla \cdot (p \nabla \frac{\delta E}{\delta p}) + \beta^{-1}\Delta p, 
\]
where \(\frac{\delta E}{\delta p}\) is the Frechet derivative of \(E\) with respect to \(p\). The particle 
dynamics is 
\[
dX = - (\nabla \frac{\delta E}{\delta p})(X) + \sqrt{2\beta^{-1}}dB_t
\]
with the drift term now depends on $p$. Our approach can be readily extended to such nonlinear Fokker-Planck equations if one simulates multiple trajectories jointly with an ensemble method (for detail, see \cite{liu2001monte}).

% \subsection*{Main algo}
% Three components:

% 1. How to collect samples. Solve SDE. Euler-Maruyama. Maybe higher-order methods should also be considered.

% 2. Project to functional tensor form using sketching.

% 3. Ready for sampling or computing observables.

% \subsection*{Extensions}

% 1. No need for a gradient structure. Any dissipative conservation law is fine.
% \[
% dX = b(X)dt + dB
% \]

% 2. We can also work with nonlinear FK-type equations. 
% The energy is 
% \[
% E(p) = U(p) + \int p \ln p
% \]
% for some nonlinear potential $U(p)$. The PDE and SDE are
% \[
% dp/dt = div (p \nabla \delta U/\delta p) + \Delta p, \quad dX = -\nabla \nabla \delta U/\delta p + dB.
% \]
% The only change is that all the particles are traced simultaneously to include the mean field nonlinear effect. 

%========================================
\section{Application to Ginzburg-Landau model}\label{sec: numerics}
% {\color{red} Not essential to the manuscript writing. Will finalize, but later}

Our main numerical treatment is dedicated to a discretized Ginzburg-Landau model in varying physical dimensions. The high dimensionality here comes from the discretization of an infinite-dimensional functional equation \cite{rosen1971functional}. For the sake of simplicity, we apply functional hierarchical tensor density estimation only at time \(T = 1\). Generally, one runs the algorithm at multiple time slots to solve for the solution across \([0, T]\). 

\paragraph{The choice of basis function}

For the Ginzburg-Landau model, the $\int_{a} |1 - (x(a))^2| da$ term ensures that it is highly unlikely for \(\lVert x \rVert_{\infty}\) to be significantly bigger than \(1\). For the discretized model, each variable $x_j$ is essentially bounded within an interval \([-B,B]\) with $B=2.5$. As a result, the samples can be effectively seen as a distribution compactly supported in the domain \([-B,B]^d\). Therefore, we can use the Fourier representation choosing a maximal degree parameter \(q\) and picking the first \(n = 2q + 1\) Fourier basis functions in the sine-cosine form in $[-B,B]$, denoted \(\{\psi_{i}\}_{i=-q}^{q}\). In effect, if one denotes \(p^{\star}\) as the target density at time \(t\in [0, T]\), then our target is to approximate \(\mathcal{P}_{q}p^{\star}\), which denotes the truncated Fourier series approximation of \(p^{\star}\) restricted to degree \(q\). 
This relationship is expressed as follows:
\begin{equation}\label{eqn: spectral representation}   
    (\mathcal{P}_{q}p^{\star})(x_1,\ldots,x_d) = \sum_{i_{1}, \ldots, i_{d} = -q}^{q} C^{\star}_{i_1,\ldots, i_d} \psi_{i_1}(x_1)\cdots \psi_{i_d}(x_d),
\end{equation}
where \(C^{\star} \in \mathbb{R}^{n^d}\) is the coefficient tensor. When the diffusion term in the Fokker-Planck equation has a strong influence, the spectral approximation in \eqref{eqn: spectral representation} provides an accurate approximation with a small \(q\). The goal is then to approximate \(p^{\star}\) through a compressed hierarchical tensor \(C \in \mathbb{R}^{n^d}\):
\begin{equation}
    p^{\star} \approx \mathcal{P}_{q}p^{\star} \approx p_{\mathrm{FHT}},\quad
    p_{\mathrm{FHT}}(x_1,\ldots,x_d)
    = \sum_{i_{1}, \ldots, i_{d} = -q}^{q} C_{i_1,\ldots, i_d} \psi_{i_1}(x_1)\cdots \psi_{i_d}(x_d).
\end{equation}

\paragraph{Bipartition structure}
Similar to the proposed hierarchical bipartition for the 2D case, we recursively bipartitions the Cartesian grid along its axes in an alternating fashion. The structure for the two-dimensional case is shown in Figure \ref{fig:2D_GZ_bipar} and can be done likewise to other dimensions. Explicitly, assume for simplicity that a function \(u\) in $\Delta$-dimensions is discretized to the \(d=m^\Delta\) points denoted \(\{u_{i_1, \ldots, i_\Delta}\}_{i_1, \ldots, i_\Delta \in [m]}\) with \(m = 2^{\mu}\). For each index \(i_\delta\), one performs a length \(\mu\) binary expansion $i_\delta = a_{\delta 1}\ldots a_{\delta \mu}$. The variable \(u_{i_1, \ldots, i_\Delta}\) is given index \(k\) with the length \(\mu \Delta\) binary expansion \(k = a_{11} a_{21} \ldots a_{\Delta-1,\mu} a_{\Delta,\mu}\), and the associated binary decomposition structure follows from the given index according to \eqref{eqn: bipartition}. One can check that the construction exactly corresponds to performing binary partition by sweeping along the $\Delta$ dimensions. For intuition on this construction, one can check that the bipartition in Figure \ref{fig:2D_GZ_bipar} corresponds to the indexing \(k = a_{11} a_{21} a_{12} a_{22} a_{13} a_{23}\).

\paragraph{The choice of sketch function}

In practice, one does not have ready access to a pre-determined sketch function, and a higher quality in the sketch function would mean that one would be better able to capture the spatial correlation of the density function. As can be seen in \eqref{eqn: Z_ag approximation}, one needs to choose relatively smooth sketch functions to ensure that an accurate approximation of terms is possible with sample estimation. Moreover, the sketch functions should be designed to capture important features so that terms such as \(Z_{a; \Bar{a}}\) could provide meaningful linear systems to solve for \(G_{q}\).

In our work, the sketch functions we choose are of two kinds. The first kind of sketch function comprises monomials of the Fourier basis of low degree. For \(a = I_{k}^{(l)}\) or \(a = [d] - I_{k}^{(l)}\) (see definition in \eqref{eqn: bipartition}), the chosen sketch functions are of the following form:
\[
f(x_a) = \prod_{j \in a}{\psi_{i_j}(x_j)},
\]
for some choice of $-q \le i_j \le q$ for each $j$. The total Fourier degree can be characterized by \(\mathrm{deg}(f) = \sum_{j\in a} |i_j| \), and one typically choose \(f\) with small \(\mathrm{deg}(f)\). 

The second kind of sketch function is motivated by the renormalization group, in particular by the coarse-graining of the variables. In particular, one chooses a cluster \(h \subset a\) and a specific Fourier mode index \(i\). The resultant function \(g_{h,i}\) is an averaging of the \(i\)-th Fourier mode over variables in \(h\):
\[
g_{h,i}(x_a) \equiv \frac{1}{|h|}\sum_{j \in h}{\psi_{i}(x_j)}.
\]
In our particular implementation, for \(a = I_{k}^{(l)},\) a typical example of chosen sketch function would be \(h = I_{2k}^{(l+1)}\) and \(i = 1\). Then, the resultant sketch function would be
\[
g_{h,1}(x_a) \equiv \frac{1}{|h|}\sum_{j \in h}{\sin(\pi x_j/B)},
\]
which can be thought of as coarse-grained information on half of the variables in \(a\). One can similarly define \(g_{a-h,1}(x_a)\). One can augment the sketch function set by including terms such as \(g_{h,1}^2\) and \(g_{a-h,1}^2\). Empirically, including such coarse-grained low-order Fourier modes greatly increases model performance and stability. 

% We choose the sketch function following one of the two criteria: (a) a Fourier series on boundary variables, and (b) an average of \(\)

\subsection{1D Ginzburg-Landau potential}

The first numerical example is a 1D Ginzburg-Landau model. The potential energy is defined as
\begin{equation}\label{eqn: 1D GZ model}
V(x_1, \ldots ,x_m) := \frac{\lambda}{2} \sum_{i=1}^{m+1}\left(\frac{x_{i} - x_{i - 1}}{h}\right)^2 + \frac{1}{4\lambda} \sum_{i = 1}^{m} \left(1 - x_i^2\right)^2,
\end{equation}
where \( h = \frac{1}{m+1} \) and \(x_0=x_{m+1}=0\). In particular, we fix \( m = 256 \), \( \lambda = 0.01 \) and \( \beta = \frac{1}{8} \). The initial condition is \(p_{0}(x) = \delta(x - (0, \ldots, 0))\) and the dimension $d$ is $256$. 

\begin{figure}[h!]
  \centering
  \includegraphics[width = 0.8\textwidth]{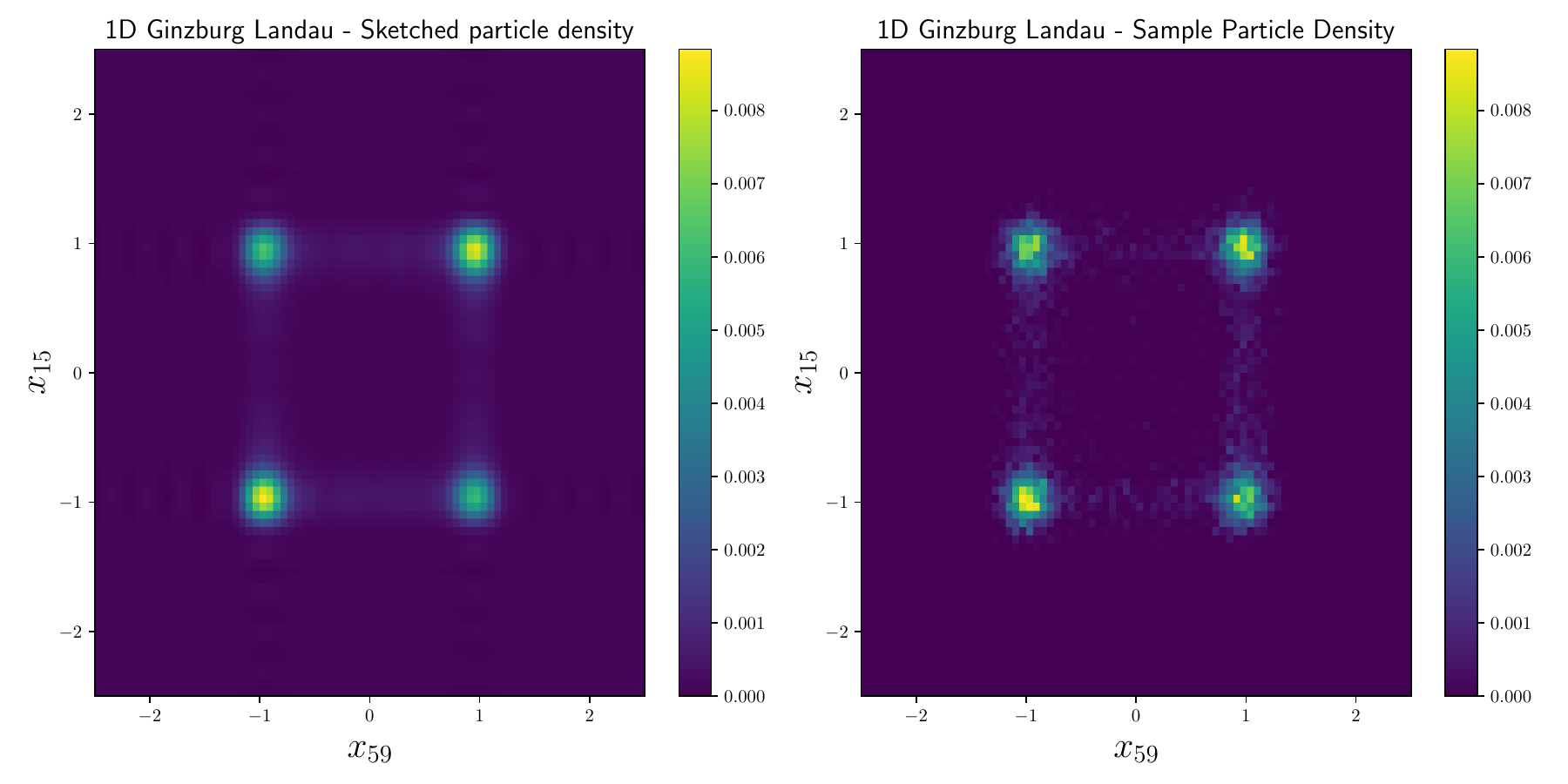}
  \caption{1D Ginzburg-Landau model. The plots of the marginal distribution at \((x_{19}, x_{59})\).}
  \label{Fig: 1D marginal}
\end{figure}

We perform \(N = 6000\) SDE simulations with \(T = 1\) and \( \delta t = \frac{1}{2000} \), all starting from the initial location $\left(0, \ldots, 0\right)$. The maximal Fourier degree \(q\) is taken to be \(q = 15\) to account for the localization of the samples, and the maximal internal bond dimension is \(r = 6\).

Figure \ref{Fig: 1D marginal} compares the marginal distribution of \((x_{15}, x_{59})\) between empirical distributions and the model obtained by hierarchical tensor sketching. One can see that hierarchical tensor sketching is successful at capturing the correlation between two variables at relatively faraway points.

\subsection{2D Ginzburg-Landau potential}
In the second numerical example, we consider the usual 2D Ginzburg-Landau model. 
The potential energy is defined as
\begin{equation}
  V(x_{(1,1)}, \ldots ,x_{(m,m)}) := \frac{\lambda}{2} \sum_{v \sim w}\left(\frac{x_{v} - x_{w}}{h}\right)^2 +  \frac{1}{4\lambda} \sum_{v}\left(1 - x_v^2\right)^2,
\end{equation}
where \(h = \frac{1}{m+1}\) and \( x_{(l, 0)} = x_{(l, m + 1)} = x_{(0, l)} = x_{(m+1, l)} = 0\) for \(l = 1, \ldots, m\). The parameters are \(m=16\), \(\lambda = 0.03\), \(\beta = \frac{1}{5}\), and the initial condition is \(p_{0}(x) = \delta(x - (0, \ldots, 0))\). The dimension $d$ is $256$. 
\begin{figure}[h!]
    \centering
    \includegraphics[width = 0.8\textwidth]{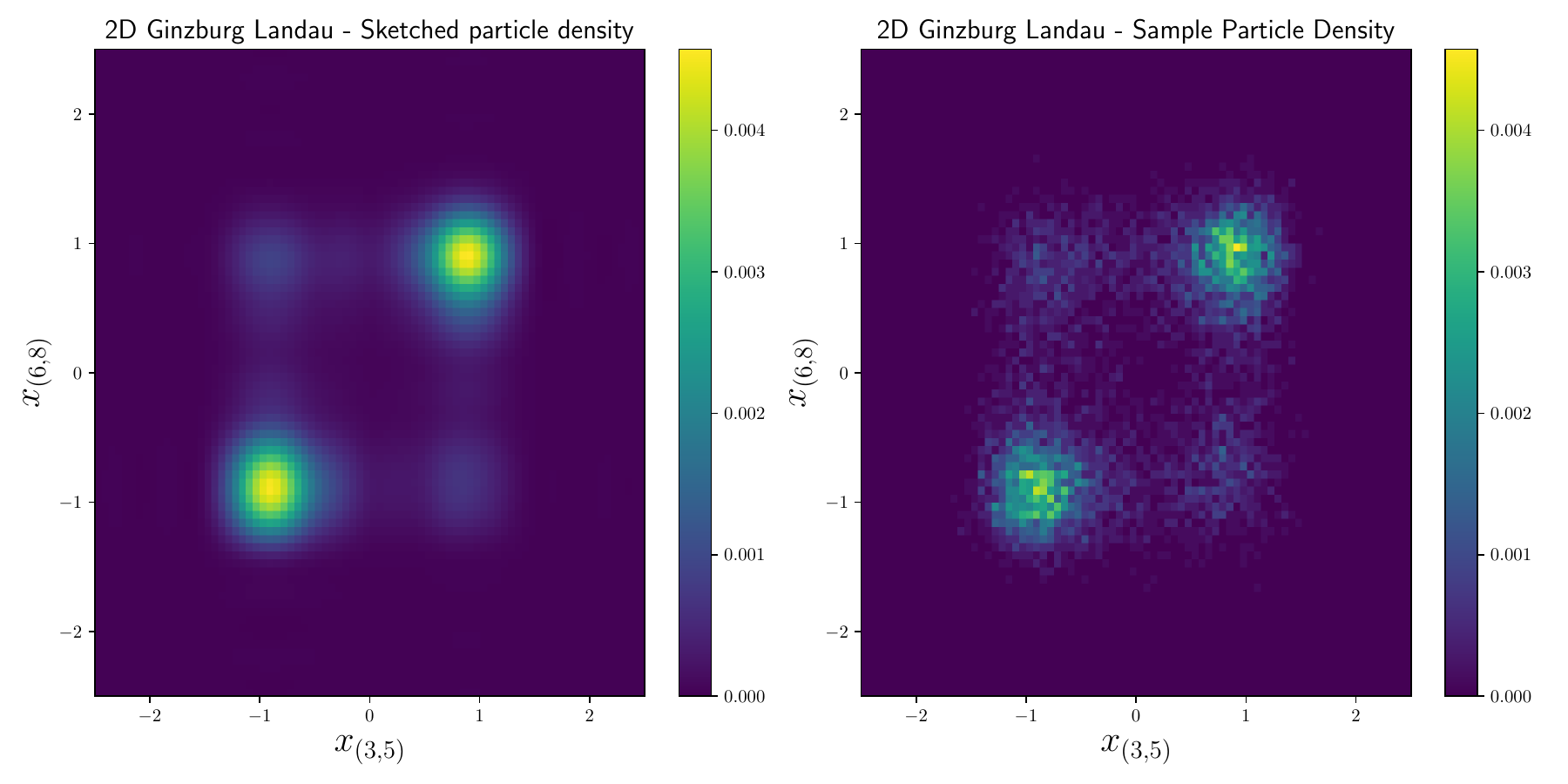}
    \caption{2D Ginzburg-Landau model. The plots of marginal distribution at \((x_{(3,5)}, x_{(6,8)})\).} 
    \label{Fig: 2D marginal}
\end{figure}

\begin{figure}[h!]
    \centering
    \includegraphics[width = 0.8\textwidth]{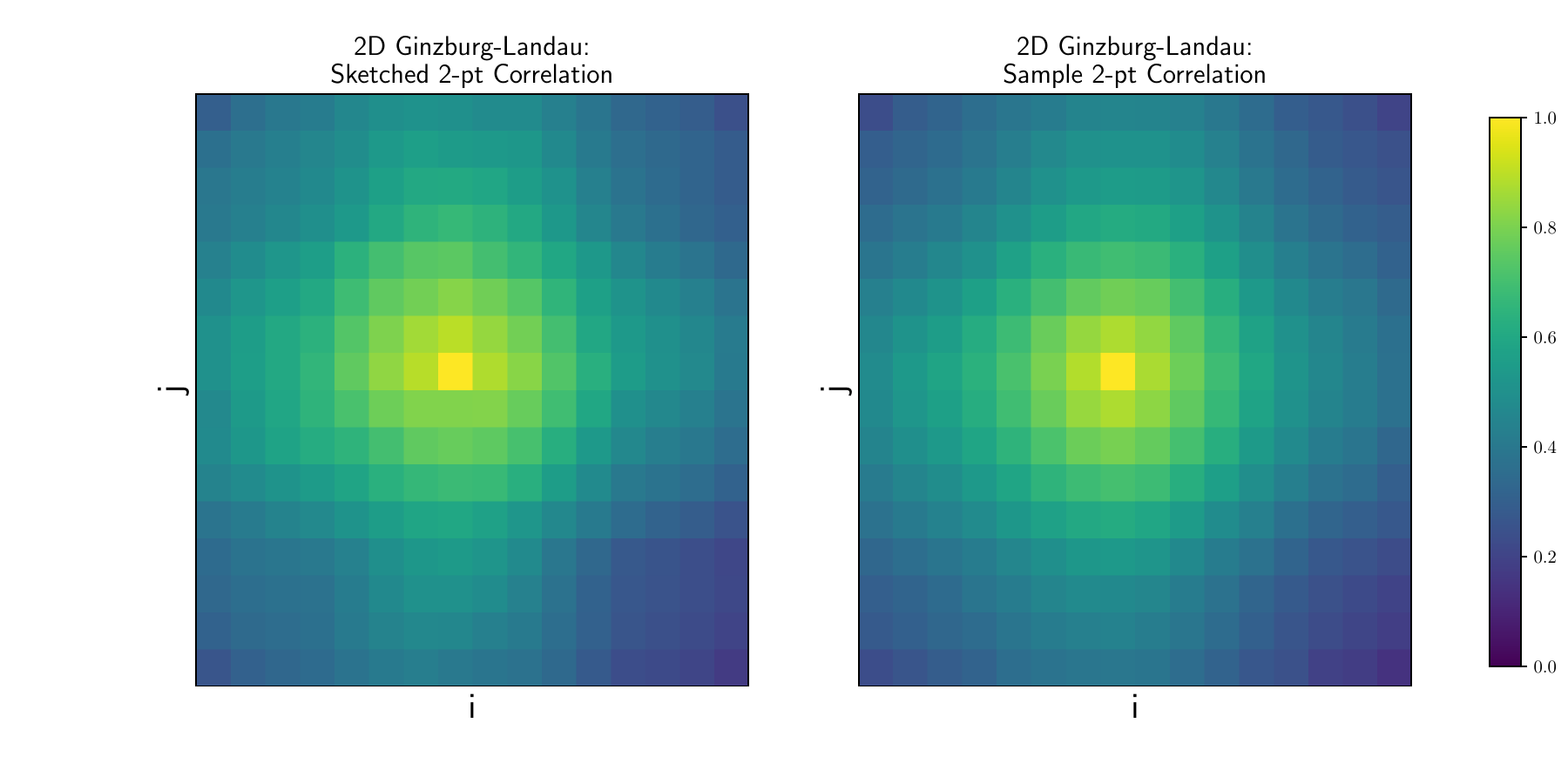}
    \caption{Two-point correlation function \(f(i,j) = \mathrm{Corr}\left(x_{(8,8)}, x_{(i,j)}\right)\) of the sketched functional hierarchical tensor, compared with the ground truth from a Monte-Carlo estimation from 60,000 samples.}
    \label{Fig: 2D_GZ_2_pt_correlation_256}
\end{figure}

We perform \(N = 6000\) SDE simulation with \(T = 1\) and \( \delta t = 0.003 \), starting from the initial condition \(\left(0, \ldots, 0\right)\). We take \(q = 10\) and a maximal internal bond of \(r = 20\). To ensure numerical stability, the internal bond dimension is dynamically chosen according to the singular values information obtained during sketching.

Figure \ref{Fig: 2D marginal} shows that the obtained marginal distribution of \((x_{(3,5)}, x_{(6, 8)})\) closely matches the empirical distributions. To test the accuracy of observable estimation, Figure \ref{Fig: 2D_GZ_2_pt_correlation_256} plots the predicted two-point correlation function 
\[
f(i, j) := \mathrm{Corr}\left(x_{(8,8)}, x_{(i,j)}\right) \equiv \frac{\mathrm{Cov}\left[x_{(8,8)}, x_{(i,j)}\right]}{\sigma_{x_{(8,8)}} \sigma_{x_{(i,j)}}},
\]
where the result closely matches the ground truth obtained from a Monte-Carlo estimation from 60,000 samples. The accuracy is quite remarkable as the internal bond \(r\) is kept at a relatively small level.

\subsection{3D Ginzburg-Landau potential}
Here, we consider a 3D Ginzburg-Landau model. Similarly, define a 3D Cartesian grid \(D = (Q, E)\), where \(Q:= \left\{(i, j,k) \mid i,j,k = 0, \ldots, m+1\right\}\). The potential energy is defined as
\begin{equation}\label{eqn: 3D GZ model}
  V(x_{(1,1,1)}, \ldots ,x_{(m,m,m)}) := \frac{\lambda}{2} \sum_{v \sim w}\left(\frac{x_{v} - x_{w}}{h}\right)^2 +  \frac{1}{4\lambda} \sum_{v}\left(1 - x_v^2\right)^2,
\end{equation}
where \(h = \frac{1}{m + 1}\), and the boundary condition is such that \(x_{(i,j,k)} = 0\) if one of \(i,j,k\) equals to \(0\) or \(m + 1\). For parameters, we fix \(m=8\), \(\lambda = 0.01\) and \(\beta = \frac{1}{10}\). The initial condition is \(p_{0}(x) = \delta(x - (0, \ldots, 0))\) and the dimension $d$ is $512$. 

\begin{figure}[h!]
    \centering
    \includegraphics[width = 0.8\textwidth]{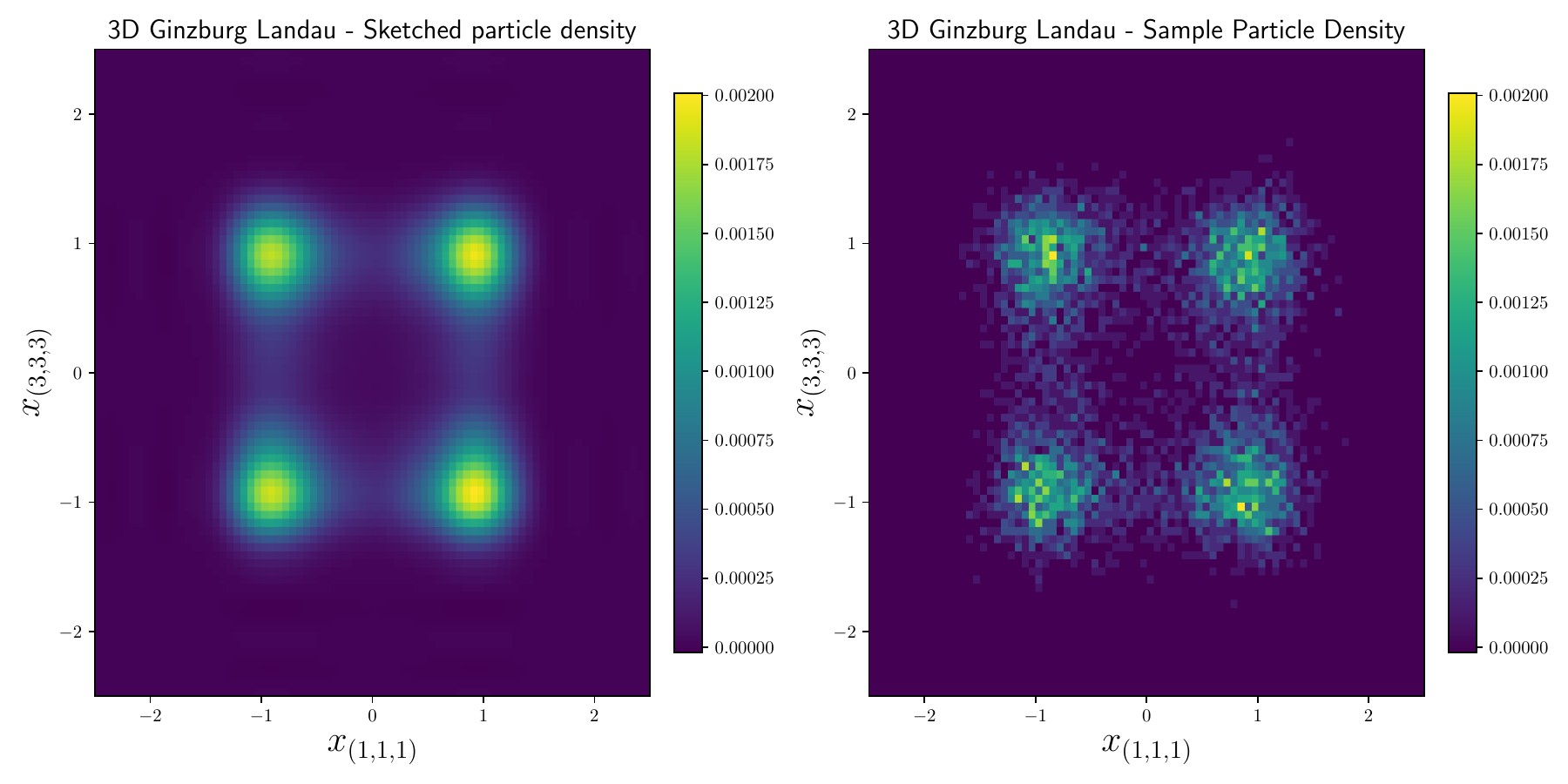}
    \caption{3D Ginzburg-Landau model. The plots of marginal distribution at \((x_{(1,1,1)}, x_{(3,3,3)})\).}     \label{Fig: 3D marginal}
\end{figure}

We perform \(N = 6000\) SDE simulation with \(T = 1\) and \( \delta t = 0.002 \), starting from the initial condition \(\left(0, \ldots, 0\right)\).  We take \(q = 10\) and a maximal internal bond of \(r = 20\), where the internal bond is likewise chosen dynamically. Figure \ref{Fig: 3D marginal} shows the performance on the marginal distribution of \((x_{(1,1,1)}, x_{(3,3,3)})\) between empirical distributions and the model obtained by hierarchical tensor sketching, and similarly the result shows quite good performance.

% In Fig \ref{}, we visualize 

%========================================
\section{Conclusion}
We introduce a novel density estimation approach that combines particle methods with a functional hierarchical tensor network in solving the high-dimensional Fokker-Planck equations. The algorithm is applied to the discretized Ginzburg-Landau model in 1D, 2D, and 3D. This method points to a new direction for tackling general high-dimensional Fokker-Planck equations. Combined with sensible numerical treatment, this approach has the potential for modeling high-dimensional particle dynamics with high fidelity. An open question is whether one can extend a tensor network approach to solving similar equations, such as the forward Kolmogorov equation and the Hamilton-Jacobi equation.

\paragraph{Acknowledgement}
We thank Yuehaw Khoo for constructive discussions.

\bibliographystyle{elsarticle-num-names} 
\bibliography{reference}

\end{document}